\documentclass[12pt,noamsfonts]{amsart}

\usepackage{amsmath}
\usepackage{amsfonts}
\usepackage{amssymb}
\usepackage{amscd}
\usepackage{verbatim}

\newcommand{\marginlabel}[1]%
  {\mbox{}\marginpar{\raggedleft\hspace{0pt}\bfseries\sf#1}}

\def\ZZ{{\mathbb Z}}
\def\NN{{\mathbb N}}

\def\CC{{\mathbb C}}
\def\AA{{\mathbb A}}
\def\RR{{\mathbb R}}
\def\QQ{{\mathbb Q}}

\def\cI{\mathcal{I}}

\def\cO{\mathcal{O}}

\DeclareMathOperator{\Hom}{Hom}

\DeclareMathOperator{\Spec}{Spec}

\DeclareMathOperator{\totdiscrep}{totdiscrep}

\newtheorem{lemma}{Lemma}[section]
\newtheorem{theorem}[lemma]{Theorem}
\newtheorem{corollary}[lemma]{Corollary}
\newtheorem{proposition}[lemma]{Proposition}

\theoremstyle{definition}
\newtheorem{definition}[lemma]{Definition}
\newtheorem{example}[lemma]{Example}

\newtheorem{remark}[lemma]{Remark}

\numberwithin{equation}{section}

\newcommand{\bean}{\begin{eqnarray}}
\newcommand{\eean}{\end{eqnarray}}
\newcommand{\be}{\begin{displaymath}}
\newcommand{\ee}{\end{displaymath}}
\newcommand{\bea}{\begin{eqnarray*}}   
\newcommand{\eea}{\end{eqnarray*}}

\begin{document}

\title{Singularities of pairs via jet schemes} 

\author[M. Musta\c{t}\v{a}]{Mircea~Musta\c{t}\v{a}}
\address{Department of Mathematics, University of California,
Berkeley, CA, 94720 and Institute of Mathematics of the
Romanian Academy}
\email{{\tt mustata@math.berkeley.edu}}

\date{\today}
\maketitle

\bigskip

\section*{Introduction}

Let $X$ be a smooth complex variety and $Y$ a closed subscheme of $X$.
The study of the singularities of the pair $(X,Y)$ is a topic which received
a lot of atention, due mainly to applications to the classification of 
higher dimensional algebraic varieties. Our main goal in this paper is
to propose a new point of view in this study, based on the properties of
the jet schemes of $Y$.

For an arbitrary scheme $W$, the $m$th jet scheme $W_m$ parametrizes morphisms
$\Spec\,\CC[t]/(t^{m+1})\longrightarrow W$. Our main result is the
 following theorem.

\begin{theorem}\label{main_result}
If $X$ is a smooth variety, $Y\subset X$ a closed subscheme, and $q>0$
a rational number, then:
\begin{enumerate}
\item The pair $(X, q\cdot Y)$ is log canonical if and only if 
$\dim\,Y_m\leq (m+1)(\dim\,X-q)$, for all $m$.
\item The pair $(X, q\cdot Y)$ is Kawamata log terminal if and only if 
$\dim\,Y_m<(m+1)(\dim\,X-q)$, for all $m$.
\end{enumerate}
\end{theorem}

In fact, it is enough to check 
 the conditions in Theorem~\ref{main_result} for one value of 
$m$, depending on a log resolution of $(X,Y)$ (see Theorem~\ref{main}
for the precise statement). As a consequence of the above result, we obtain
a formula for the log canonical threshold.

\begin{corollary}\label{main_corollary}
If $X$ is a smooth variety and $Y\subset X$ is a closed subscheme, then the
log canonical threshold of the pair $(X,Y)$ is given by
$$c(X,Y)=\dim\,X-\sup_{m\geq 0}{{\dim\,Y_m}\over{m+1}}.$$
\end{corollary}

In the same spirit with the previous remark, the above supremum can be obtained
for specific values of $m$.
A similar formula holds for the log canonical threshold around a closed
subset $Z\subset X$. Note that as a consequence of the above formula, 
the log canonical threshold depends only on $Y$ and $\dim\,X$, but not
on the particular ambient variety or the embedding. 

We apply Corollary~\ref{main_corollary} to give simpler proofs of some
results on the log canonical threshold proved by Demailly and Koll\'{a}r
in \cite{demailly} using analytic techniques. For example, we use a 
semicontinuity statement about the dimension of the jet schemes to deduce
the semicontinuity of the log canonical threshold.

The main technique we use in the proof of Theorem~\ref{main_result}
is motivic integration, technique due to Kontsevich, Batyrev, and Denef
and Loeser. The idea is similar to the one used in our previous paper
\cite{mustata}. Here is a brief description of the idea of the proof.

Let $X_{\infty}=\projlim_m X_m$ be the space of arcs of $X$. A $\CC$-valued
point of $X_{\infty}$ corresponds to a morphism $\Spec\,\CC[[t]]
\longrightarrow X$. Let $F_Y\,:\,X_{\infty}\longrightarrow\NN\cup\{\infty\}$
be the function such that for an arc $\gamma$ over $x\in X$ 
which corresponds to a morphism $\widetilde{\gamma}\,:\,{\mathcal O}_{X,x}
\longrightarrow\CC[[t]]$, we have
$F_Y(\gamma)={\rm ord}\,(\widetilde{\gamma}({\mathcal I}_{Y,x}))$.
For a suitable function $f\,:\,\NN
\longrightarrow\NN$, the motivic integral of $f\circ F_Y$ on 
$X_{\infty}$ encodes the information about the
dimension of $Y_m$, for all $m$.

On the other hand, if we take a log resolution $\pi\,:\,X'
\longrightarrow X$, then by the change of variable formula of Batyrev
(\cite{batyrev}) and Denef and Loeser (\cite{denef}), this integral can be
 expressed as the integral on $X'_{\infty}$ of $f\circ F_{\pi^{-1}(Y)}
+F_W$. Here $F_{\pi^{-1}(Y)}$ and $F_W$ are the functions associated
as above with the effective divisors $\pi^{-1}(Y)$ and $W$, where $W$ is the
relative canonical divisor of $\pi$. Since both $W$ and $\pi^{-1}(Y)$
are divisors with simple normal crossings, this last integral can be
 explicitely computed by a formula involving the coefficients of $W$ and
$\pi^{-1}(Y)$. By an inspection of the monomials which enter in the expression
of the two integrals we deduce the equivalence in Theorem~\ref{main_result}.

\smallskip

The paper is organized as follows. In the first section we review the
 definitions of log canonical and Kawamata log terminal pairs.
  In the
 algebraic context the usual setting in the literature 
 is that of a pair $(X,Y)$, where
$Y$ is a divisor on $X$. 
Therefore we give the extension to the case of pairs of arbitrary
codimension, making reference for the case of divisors to Koll\'{a}r
(\cite{kollar}). Note, however, that we will always stick to the
case of a smooth ambient variety, as our main results are valid
only in this context.

The second section contains an overview of the basic definitions and 
notions related to jet schemes and motivic integration. We give also some 
further properties of jet schemes. In particular, we prove a 
semicontinuity result about the dimensions of the jet schemes.

The third section is the heart of the paper. It gives the proof of
Theorem~\ref{main_result} along the lines described above. In the last section
we show how Corollary~\ref{main_corollary} can be used to prove properties
of the log canonical threshold. For example, the canonical isomorphism
$(Y\times Y')_m\simeq Y_m\times Y'_m$ implies that
$c(X\times X',Y\times Y')=c(X,Y)+c(X'.Y')$. This can be used to prove 
the inequality in \cite{demailly}: $c_z(X,Y_1\cap Y_2)\leq
c_z(X,Y_1)+c_z(X,Y_2)$, for every $z\in Y_1\cap Y_2$. Another application
is to the semicontinuity theorem in \cite{demailly} for
log canonical thresholds. We conclude with a concrete example: we use
the formula in Corollary~\ref{main_corollary} to compute the log
canonical threshold of monomial ideals, recovering in this way a result
from \cite{agv} (see also \cite{howald}
for a proof via multiplier ideals).

\smallskip

\subsection{Acknowledgements}
I am grateful to David Eisenbud and Edward Frenkel for introducing me
to jet schemes and for their generous help and support. I am indebted to
Lawrence Ein for pointing out an incomplete argument in a previous version
of this paper.

\bigskip
\bigskip

\section{Log terminal and log canonical pairs}

\bigskip

All schemes are of finite type over $\CC$. A variety is a reduced,
irreducible scheme.
We use the theory of singularities of pairs for which our main reference
is \cite{kollar}. However, as we work with pairs of arbitrary codimension,
we review in this section some extensions to this setting of the definitions
we need. For an analytic approach to singularities of pairs of
arbitrary codimension, see \cite{demailly}.

We consider pairs $(X,Y)$, where
$X$ is a smooth variety and  $Y\hookrightarrow X$ is a closed
subscheme with $Y\neq X$. 
 A log resolution of $(X, Y)$ is a proper, birational
morphism $\pi\,:\,X'\longrightarrow X$ such that $X'$ is smooth,
 $\pi^{-1}(Y)=D$ is an effective divisor and the union
$D\cup Ex(\pi)$ has simple normal crossings. Here $Ex(\pi)$
denotes the exceptional locus of $\pi$.
The results in \cite{hironaka} show that log resolutions exist
(in fact, one can assume in addition that $\pi$ is an isomorphism
over $X\setminus Y$.)

For an arbitrary proper, birational morphism
$\pi\,:\,X'\longrightarrow X$, with $X'$ and $X$ smooth,
the relative canonical divisor $K_{X'/X}$ is the unique effective divisor
supported on $Ex(\pi)$ such that $\cO_{X'}(K_{X'/X})\simeq
\omega_{X'}\otimes\pi^*(\omega_X^{-1})$.

\begin{definition}
Let $(X,Y)$ be a pair as above and $q>0$ a rational number.
For a log resolution $\pi\,:\,X'\longrightarrow X$ of $(X,Y)$, 
with $\pi^{-1}(Y)=D$, we write $K_{X'/X}-q\cdot D
=\sum_i\alpha_iE_i$. The total discrepancy of
 $(X, q\cdot Y)$ is defined by
\begin{equation}
{\rm totdiscrep}\,(X,q\cdot Y)=
\begin{cases}
-\infty, &\text{if $\alpha_i<-1$, for some $i$;}\\
\min_i\{0,\alpha_i\}, &\text{otherwise.}\\
\end{cases}
\end{equation}
The pair $(X, q\cdot Y)$ is called Kawamata log terminal or log canonical if
${\rm totdiscrep}\,(X, q\cdot Y)>-1$ or
${\rm totdiscrep}\,(X,q\cdot Y)\geq -1$, respectively.
\end{definition}

\begin{remark}\label{connection}
Since $q\cdot D-K_{X'/X}$ is a divisor with simple normal crossings,
it follows from Lemma 3.11 in \cite{kollar}, that we could have defined
the total discrepancy of $(X,q\cdot Y)$ as the total discrepancy
of $q\cdot D-K_{X'/X}$. In particular, $(X,q\cdot Y)$ is
Kawamata log terminal
(log canonical) if and only if $(X', q\cdot D-K_{X'/X})$ is Kawamata log
terminal (log canonical). 
\end{remark}

\begin{proposition}\label{no_dependance}
The total discrepancy does not depend on the particular log resolution,
and therefore neither do the notions of 
log canonical and Kawamata log terminal
pairs.
\end{proposition}

\begin{proof}
By the above remark,
it is enough to show that if ${\pi'}\,:\,X''\longrightarrow X'$ is a log
resolution of $(X',D\cup Ex(\pi))$, with ${\pi'}^{-1}(D)=F$, then 
$$\totdiscrep\,(q\cdot D-K_{X'/X})=\totdiscrep\,(q\cdot F-K_{X''/X}).$$
This follows from the fact that
$K_{X''/X}=K_{X''/X'}+{\pi'}^*(K_{X'/X})$
and the invariance of the total discrepancy under proper birational
morphisms in the case of divisors (see \cite{kollar}, Lemma 3.10).
\end{proof}

\begin{definition}
Let $(X, Y)$ be a pair as above, $Z\subset X$ a nonempty closed subset,
 and $q>0$ a rational number. The log canonical threshold of $(X, q\cdot Y)$
around $Z$ is defined by $c_Z(X, q\cdot Y)=\sup\{q'\in\QQ_+\mid
(X, qq'\cdot Y)\,\text{is log canonical in an open neighbourhood of}\,Z\}$.
When $Z=X$, we omit it from the notation.
\end{definition}

As it is obvious that $c_Z(X, q\cdot Y)=(1/q)\,c(X, Y)$, it follows that
when computing log canonical thresholds,
there is no loss of generality in assuming that 
$q=1$. The proposition below is the analogue 
of Proposition 8.5   in \cite{kollar}. It shows
 how to compute the log canonical threshold
from a log resolution.

Let $(X, Y)$ be a pair and  $Z\subseteq X$ a nonempty closed subset. Let
$\pi\,:\,X'\longrightarrow X$ be a log resolution of $(X, Y)$ which is
an isomorphism over $X\setminus Y$. We write $\pi^{-1}(Y)=D=
\sum_ia_iE_i$, with $a_i\geq 1$ and $K_{X'/X}=\sum_ib_iE_i$.

\begin{proposition}(\,\cite{demailly} 1.7, \cite{kollar} 8.5)\label{lc_formula}
With the above notation, the log canonical threshold of $(X,Y)$
around $Z$ is given by
\begin{equation}
c_Z(X,Y)=\min_i\{(b_i+1)/a_i\,\mid\,\pi(E_i)\cap Z\neq\emptyset\}.
\end{equation}
In particular, either $c_Z(X,Y)$ is a positive rational number,
or $\,\,\,\,$
 $c_Z(X,Y)=\infty$.
\end{proposition}

\begin{proof}
We have $(X, q\cdot Y)$ log canonical if and only if 
$b_i-qa_i\geq -1$, 
or equivalently $q\leq (b_i+1)/a_i$,
for all $i$. By replacing $X$ with 
the complement of the union of those images of $E_i$ which do not intersect 
$Z$, we get our result.
\end{proof}

\begin{remark}
The log canonical threshold $c_Z(X,Y)$ is infinite if and only if 
$Z\cap Y=\emptyset$.
\end{remark}

\begin{remark}
We assumed above that for every pair $(X,Y)$, we have $Y\neq X$.
We will sometimes find convenient to include also the case of the pair 
$(X,X)$, so that we make the convention $c_Z(X,X)=0$, for every nonempty
closed subset $Z\subseteq X$.
\end{remark}

Though unnecessary for our needs, for the sake of completeness
we mention the translation of the above notions into the language
of multiplier ideals (see \cite{ein} or \cite{lazarsfeld}). Recall that if 
$\pi$ is a log resolution of $(X, Y)$ such that $\pi^{-1}(Y)=D$, then
the multiplier ideal of the pair
$(X,q\cdot Y)$ is defined by
\begin{equation}
I(X, q\cdot\cI_{Y/X})=\pi_*(\cO_{X'}(K_{X'/X}-[q\cdot D])),
\end{equation}
where $[q\cdot D]$ denotes the integral part of the $\QQ$--divisor
$q\cdot D$.

\begin{proposition}
The pair $(X,q\cdot Y)$ is Kawamata log terminal if and only if 
$I(X, q\cdot\cI_{Y/X})=\cO_X$ and it is log canonical if and only if
$I(X, q'\cdot \cI_{Y/X})=\cO_X$ for every $q'<q$ . In particular,
if $Z\subseteq X$ is a nonempty closed subset, then
\begin{equation}
c\,(X,Y)=\sup\{q>0 \mid I(X,q\cdot\cI_{Y/X})=\cO_X\,{\rm around}\,Z\}.
\end{equation}
\end{proposition}

\begin{proof}
It is enough to prove the first assertion. This follows from the fact that
$I(X,q\cdot\cI_{Y/X})=\cO_X$ if and only if $K_{X'/X}-[q\cdot D]$
is effective or equivalently, all the coefficients in
$K_{X'/X}-q\cdot D$ are $>-1$.
\end{proof}

\bigskip
\bigskip

\section{An overview of jet schemes and motivic integration}

\bigskip

In this section we collect the definitions and basic properties of jet
schemes. We review also motivic integration on smooth varieties, technique
which plays a major role in the proof of our main result in the next section.

 For every scheme  $W$
(of finite type over $\CC$) and every $m\in\NN$, the jet scheme $W_m$
is a scheme of finite type over $\CC$ characterized by
\begin{equation}
\Hom(\Spec\,A,W_m)\simeq\Hom(\Spec\,A[t]/(t^{m+1}), W),
\end{equation}
for every $\CC$--algebra $A$. 
 In particular, the closed points
of $W_m$ correspond to $\CC[t]/(t^{m+1})$--valued points of $W$.
The correspondence $W\longrightarrow W_m$ is a functor. In fact, from the
definition it follows that it is the right adjoint of the functor
$Z\longrightarrow Z\times\Spec\,\CC[t]/(t^{m+1})$.

It is clear that we have $W_0\simeq W$ and $W_1\simeq T\,W$, the tangent
space of $W$. In general, we have projections $\phi_m\,:\,W_m\longrightarrow
W_{m-1}$ which are induced by the projections $A[t]/(t^{m+1})
\longrightarrow A[t]/(t^m)$, for a $\CC$--algebra $A$.
By composing these morphisms we get natural projections
$\rho_m\,:\,W_m\longrightarrow W$. Whenever the variety is not clear from the
context, we will add a superscript: for example, $\rho_m^W$.
The projective limit of the schemes $W_m$ is a scheme 
$W_{\infty}$, in general not of finite type over $\CC$,
called the space of arcs of $W$. Its $\CC$--valued points
correspond to $\CC[[t]]$--valued points of $W$. We have canonical
morphisms $\psi_m\,:\,W_{\infty}\longrightarrow W_m$.

This construction is local in the sense that if $U\subseteq W$ is an
open subset, then $U_m\simeq \rho_m^{-1}(U)$, for all $m$. 
More generally, if $f\,:\,W'\longrightarrow W$ is an \'{e}tale
morphism, then $W'_m\simeq W_m\times_WW'$.

 We can therefore
reduce the description of $W_m$ to the case when $W$ is affine.
In this case one can write down explicit equations
for $W_m$ as follows. Suppose that
$W\subseteq\Spec\,\CC[X_i;i\in I]$ is defined by a system of polynomials
$(f_{\alpha})_{\alpha}$. Consider the ring $R_m=\CC[X_i,X'_i,\ldots X_i^{(m)};
i\in I]$ and $D$ the unique $\CC$--derivation of $R_m$ such that 
$D(X^{(j)}_i)=X^{(j+1)}_i$ for all $i$ and $j$, where $X_i^{(0)}=X_i$
and $X_i^{(m+1)}=0$. For a polynomial $f\in \CC[X_i; i\in I]$, we put
$f^{(j)}=D^j(f)$. With this notation, $W_m\subseteq\Spec\,R_m$
is defined by $(f_{\alpha},f'_{\alpha},\ldots,f^{(m)}_{\alpha})_{\alpha}$. 
In follows from this description that if $W'\hookrightarrow
 W$ is a closed subscheme of $W$, then the induced morphism 
$W'_m\longrightarrow W_m$ is also a closed immersion.

If $W$ is a smooth variety of dimension $n$,
 then the morphisms $\phi_m$ are locally trivial with fiber $\AA^n$.
In particular, $W_m$ is a smooth variety of dimension 
$(m+1)n$.

\bigskip

We discuss some more properties of jet schemes which we will
apply in the last section to deduce corresponding properties
of the log canonical threshold. 

\begin{proposition}\label{product_of_jets}
For every two schemes $W$ and $Z$ and every $m$, there is
a natural isomorphism $(W\times Z)_m\simeq W_m\times Z_m$.
\end{proposition}

\begin{proof}
The assertion follows from the fact that the functor
$W\longrightarrow W_m$ has a left adjoint, and therefore
it commutes with direct products.
\end{proof}

For every scheme $W$ and every $m\geq 1$ there is a natural ``action'':
$$\Phi_m\,:\,{\bf A}^1\times W_m\longrightarrow W_m$$
over $W$ defined at the level of $A$-valued points as follows. 
For an algebra $A$, an $A$-valued point in ${\bf A}^1\times W_m$ consists
of a pair $(a,f)$, for some
 $a\in A$ and a morphism $f\,:\,\Spec\,A[t]/(t^{m+1})
\longrightarrow W$.
 $\Phi_m(a,f)$ is the composition $f\circ g_a$, where $g_a$ is 
induced by the $A$-algebra homomorphism $A[t]/(t^{m+1})
\longrightarrow A[t]/(t^{m+1})$ which maps $t$ to $at$.

It is clear that these ``actions'' are compatible with the projection
$\phi_m$. The restriction of $\Phi_m$ to $\CC^*\times W_m$ induces
an action of the torus on $W_m$. On the other hand, note that 
$\Phi_m(\{0\}\times W_m)$ is the image of the
``zero section'' $\sigma_m$ of $\rho_m$.
This is defined by composition with the scheme morphism induced
by the inclusion
$A\hookrightarrow A[t]/(t^{m+1})$. It is clear that
$W_m\setminus \sigma_m(W)$ is $\CC^*$-invariant. With this notation, we have

\begin{lemma}\label{cone}
For every $w\in W$, the fiber $\rho^{-1}_m(w)$ is the cone over the
(possibly empty) projective scheme $(\rho^{-1}_m(w)\setminus\{\sigma
_m(w)\})/{\CC^*}$.
\end{lemma}

\begin{proof}
By restricting to an open neighbourhood of $w$,
we may assume that $W\subseteq {\bf A}^N$ and
that $w=0$ is the origin. We have an embedding $W_m\subseteq {\bf A}^{(m+1)N}$
which induces an embedding $\rho_m^{-1}(0)\subseteq {\bf A}^{mN}$
such that $\sigma_m(0)$ corresponds to the origin. The action of $\CC^*$
on $\rho_m^{-1}(0)\setminus\{0\}$ extends to an action on
${\bf A}^{mN}=\Spec[X'_i,\ldots,X^{(m)}_i;i]$ induced by $\lambda\cdot
X_i^{(j)}=\lambda^jX_i^{(j)}$. Therefore $\rho_m^{-1}(0)\setminus\{0\}/
\CC^*$ is a subscheme of a weighted projective space, hence it is projective,
and $\rho_m^{-1}(w)$ is the cone over it.
\end{proof}

Our next goal is the proof of the following semicontinuity result.
For a scheme $\pi\,:\,{\mathcal W}\longrightarrow S$ over $S$
and a closed point $s\in S$, we denote the fiber of $\pi$ over $s$
by ${\mathcal W}_s$.

\begin{proposition}\label{semicontinuity1}
Let $\pi\,:\,{\mathcal W}\longrightarrow S$ be a family of schemes
and $\tau\,:\,S\longrightarrow {\mathcal W}$ a section of $\pi$.
For every $m\geq 1$, the function
$$f(s)=\dim\,(\rho_m^{{\mathcal W}_s})^{-1}(\tau(s))$$
is upper semicontinuous on the set of closed points of $S$.
\end{proposition}

\vfill\eject

The key to the proof of Proposition~\ref{semicontinuity1}
is to use a relative version of the jet schemes. Suppose we work
over a fixed scheme $S$. If ${\mathcal W}\longrightarrow S$
is a scheme over $S$, then the $m$th relative jet scheme
$({\mathcal W}/S)_m$ is characterized by 
$${\rm Hom}_S({\mathcal Z}\times\Spec\,\CC[t]/(t^{m+1}), {\mathcal W})
\simeq {\rm Hom}_S({\mathcal Z}, ({\mathcal W}/S)_m),$$
for every scheme ${\mathcal Z}$ over $S$. Therefore the functor
${\mathcal W}\longrightarrow ({\mathcal W}/S)_m$ is the right adjoint
of the functor ${\mathcal Z}\longrightarrow {\mathcal Z}\times
\Spec\,\CC[t]/(t^{m+1})$ between schemes over $S$. 

The existence of $({\mathcal W}/S)_m$ can be settled as in the 
absolute case by giving local equations. More precisely, since 
the construction is local on ${\mathcal W}$, we may assume that
both $S$ and ${\mathcal W}$ are affine: $S=\Spec\,A$ and
${\mathcal W}\hookrightarrow\Spec\,A[X_i;i]$ is defined by
$(f_{\alpha})_{\alpha}$. Then $({\mathcal W}/S)_m
\hookrightarrow\Spec\,A[X_i,X'_i,\ldots,X_i^{(m)};i]$
is defined by $(f_{\alpha}, f'_{\alpha},\ldots,f^{(m)}_{\alpha})_{\alpha}$.
Here $f^{(j)}=D^j(f)$, where $D$ is the unique derivation over $A$ of
$A[X_i,X'_i,\ldots, X^{(m)}_i;i]$ such that $D(X_i^{(p)})=
X_i^{(p+1)}$ for all $p$
(we put $X_i^{(0)}=X_i$ and $X_i^{(m+1)}=0$).

As in the absolute case, we have canonical projections
$\rho_m\,:\,({\mathcal W}/S)_m\longrightarrow {\mathcal W}$.
We have also a ``zero section'' $\sigma_m\,:\,{\mathcal W}
\longrightarrow ({\mathcal W}/S)_m$.

The following lemma follows immediately from the functorial
definition of relative jet schemes.

\begin{lemma}
For every scheme morphism $S'\longrightarrow S$,
if ${\mathcal W'}\simeq {\mathcal W}\times_SS'$, then we have a canonical
isomorphism $({\mathcal W'}/S')_m\simeq ({\mathcal W}/S)\times_SS'$,
for every $m$. In particular, for every closed point $s\in S$, the fiber
of $({\mathcal W}/S)_m$ over $s$ is isomorphic with
$({\mathcal W}_s)_m$.
\end{lemma}

Over every closed point $s\in S$, the projection $\rho_m$ and
the ``zero section'' $\sigma_m$ are the ones we defined in the
 absolute case. On the other hand, the actions of $\CC^*$
on each fiber globalize to an action on $({\mathcal W}/S)_m$
such that $({\mathcal W}/S)_m\setminus\sigma_m({\mathcal W})$
is invariant. After these preparations we can give the proof 
of Proposition~\ref{semicontinuity1}.

\begin{proof}[Proof of Proposition~\ref{semicontinuity1}]
Consider the projection $({\mathcal W}/S)_m\setminus\sigma_m({\mathcal W})
\longrightarrow {\mathcal W}$ and let ${\mathcal Z}$ be the inverse image
by this morphism of $\tau(S)$. We have a natural action of $\CC^*$ on
${\mathcal Z}$ and the quotient by this action is a scheme $\overline
{\mathcal Z}$ which is proper over $S$. The properness follows from the
fact that it is locally projective, an assertion which is just the
globalization of Lemma~\ref{cone} and can be proved in the same way. 

But for every closed point $s\in S$ we have $f(s)=\dim\,\overline{\mathcal Z}
_s+1$ (or $f(s)=0$ if $\overline{\mathcal Z}_s=\emptyset$)
 and our assertion follows from the semicontinuity of the dimension
of the fibers of a proper morphism. 
\end{proof}

\bigskip

We review now the basic facts about motivic integration on smooth varieties.
This technique was developed by 
Kontsevich \cite{kontsevich}, Batyrev \cite{batyrev} and Denef and Loeser
\cite{denef}. There are several possible extensions (see \cite{denef}
for the general treatment), but we use only
Hodge realizations of
 motivic integrals on the space of arcs of 
a smooth variety. For a nice introduction to these ideas 
we refer to Craw \cite{craw}.

 Suppose from now on that $X$ is a smooth
variety of dimension $n$. On $X_{\infty}$ the theory provides an
algebra of functions ${\mathcal M}$ and a finitely additive measure
on $\mathcal M$ with values in the Laurent power series ring in $u^{-1}$
and $v^{-1}$:
$$\mu\,:\,{\mathcal M}\longrightarrow S=\ZZ[[u^{-1},v^{-1}]][u,v].$$
On $S$ we consider the topology defined by the descending sequence of
subgroups $\{\oplus_{i+j\geq l}\ZZ u^{-i}v^{-j}\}_l$.

We will be interested in the subalgebra ${\mathcal Cyl}$ of ${\mathcal M}$
consisting of cylinders of the form: $\psi_m^{-1}(C)$, for some $m\geq 1$
and some constructible subset $C\subset X_m$. The measure of such a cylinder
is given by $\mu(\psi_m^{-1}(C))=E(C;u,v)(uv)^{-(m+1)n}$, where
$E(C;u,v)$ is the Hodge-Deligne polynomial of $C$. What is important for us
is that  if $C\subseteq X_m$ is a locally closed
subset, then $E(C; u,v)$ is a polynomial of degree $2\,(\dim\, C)$ and
the term of degree $2\,(\dim\,C)$ is $l(C)(uv)^{\dim\,C}$, where $l(C)$ is
the number of irreducible components of $C$ of maximal dimension.

Besides the cylinders, some sets of measure zero appear in the
definition of measurable functions. 
If $T\subset X_{\infty}$ is such that there is a sequence of cylinders
$W_r\in {\mathcal Cyl}$ with $T\subset W_r$ for all $r$
and $\mu(W_r)\longrightarrow 0$, then $T\in {\mathcal M}$ and
$\mu(T)=0$.

 A function $F\,:\,X_{\infty}\longrightarrow \NN\cup\{\infty\}$
 is called measurable 
if $F^{-1}(s)\in {\mathcal M}$ for every
$s\in\NN\cup\{\infty\}$ and if $\mu(F^{-1}(\infty))=0$.
If the sum $\,\,\,\,$
$\sum_{s\in\NN}\mu(F^{-1}(s)) (uv)^{-s}$ is convergent in $S$, then
$F$ is caled integrable, and the sum is called the
motivic integral of $F$ and denoted by $\int_{X_{\infty}}e^{-F}$. 

\medskip

Every subscheme $Y\hookrightarrow X$ defines a function $F_Y$
on $X_{\infty}$, as follows. If $\gamma\in X_{\infty}$ is an arc over
$x\in X$, then it can be identified with a ring homomorphism
$\widetilde{\gamma}\,:\,\cO_{X,x}\longrightarrow\CC[[t]]$. We define
$F_Y(\gamma):={\rm ord}
(\widetilde{\gamma}(\cI_{Y,x}))$. It follows from this definition
that $$F_Y^{-1}(s)=\psi_{s-1}^{-1}(Y_{s-1})\setminus
\psi_s^{-1}(Y_s),$$
for every integer $s\geq 0$ and $F_Y^{-1}(\infty)=Y_{\infty}$.
Here we make the convention $Y_{-1}=Y$ and $\psi_{-1}=\psi_0$.
It follows that $F_Y^{-1}(s)\in {\mathcal Cyl}$ for every $s\geq 0$.
Moreover, the following lemma shows that if $Y\neq X$, then 
$F_Y^{-1}(\infty)\in {\mathcal M}$ and $\mu(F_Y^{-1}(\infty))=0$.
Therefore $F_Y$ is measurable.

\begin{lemma}(\cite{mustata}, 3.7)\label{measurable}
If $D\subset X$ is a divisor and $x\in X$ is a point such that
${\rm mult}_xD=a$, then 
\begin{equation}
\dim\,(\rho_m^D)^{-1}(x)\leq nm-[m/a],
\end{equation}
where $[y]$ denotes the integral part of $y$. In particular, for every
closed subscheme $Y$ of $X$, with $Y\neq X$, we have 
$\dim\,Y_m-n(m+1)\longrightarrow -\infty$, so that $F_Y^{-1}(\infty)\in
{\mathcal M}$ and $\mu(F_Y^{-1}(\infty))=0$.
\end{lemma}

\medskip

One of the main results of the theory is the change of variable formula
for motivic integrals.

\begin{proposition}\label{change_of_variable}[\,\cite{batyrev} 6.27,
\cite{denef} 3.3]
Let $\pi\,:\,X'\longrightarrow X$ be a proper, birational morphism of
smooth varieties. Let $W=K_{X'/X}$ be the relative canonical divisor
and $\pi_{\infty}\,:\,X'_{\infty}
\longrightarrow X_{\infty}$ the morphism induced
by $\pi$. For every measurable function $F\,:\,X_{\infty}\longrightarrow
\NN\cup\{\infty\}$, we have
$$\int_{X_{\infty}}e^{-F}=\int_{X'_{\infty}}e^{-(F\circ\pi_{\infty}+F_W)},$$
meaning that one integral exists if and only if the other one does,
and in this case they are equal.
\end{proposition}

\begin{remark}
In connection with the change of variable formula, note that if $F=F_Y$,
for some proper subscheme $Y\hookrightarrow X$, then 
$F\circ\pi_{\infty}=F_{\pi^{-1}(Y)}$.
\end{remark} 

\bigskip
\bigskip

\section{Jet schemes of log terminal and log canonical pairs}

\bigskip

The main goal of this section is to prove the characterization of log canonical
and Kawamata log terminal pairs in terms of jet schemes. First let us fix the
notation.

Consider a pair $(X,Y)$ with $\dim\,X=n$ and fix a log resolution of
 this pair
$\pi\,:\,X'\longrightarrow X$ such that $\pi$ is an isomorphism over 
$X\setminus Y$. 
 Let $D=\pi^{-1}(Y)$ and
$W=K_{X'/X}$ and write
$D=\sum_{i=1}^ra_iE_i$, with $a_i\geq 1$ for all $i$
 and $W=\sum_{i=1}^rb_iE_i$.

\begin{theorem}\label{main}
With the above notation, if $q>0$,  we have the following equivalences.
\begin{enumerate}
\item $(X, q\cdot Y)$ is log canonical if and only if 
$\dim\,Y_m\leq (m+1)(n-q)$, for every $m\geq 0$.  Moreover, it is enough
to have this inequality for some $m\geq 0$ such that $a_i\mid (m+1)$
for all $i$.
\item $(X, q\cdot Y)$ is Kawamata log terminal if and only if
$\dim\,Y_m < (m+1)(n-q)$, for every $m\geq 0$. As above, it is enough
to have this inequality for some $m$ such that $a_i\mid (m+1)$ for
all $i$.
\end{enumerate}
\end{theorem}

\begin{proof}
The idea of the proof is similar to that of Theorem~3.1 in
\cite{mustata}, so that we refer to that paper for some of the details.

With the notation in the theorem, note that we have $(X, q\cdot Y)$
log canonical (Kawamata log terminal) if and only if
$b_i+1\geq (>) qa_i$,
 for all $i$. Since the case $Y=\emptyset$ is trivial
(we follow the convention that $\dim(\emptyset)=-\infty$), we assume from
now on $Y$ nonempty.

We fix a function $f\,:\,\NN\longrightarrow\NN$, such that
for every $s\geq 0$,
\begin{equation}
(\star)\,f(s+1)>f(s)+\dim\,Y_s+C(s+1),
\end{equation} where $C\in\NN$ is a constant such that
$C>\mid N-(b_i+1)/a_i\mid$, for all $i$. We extend this function by defining
$f(\infty)=\infty$. For the proof of the ``if'' part, we will add later
an extra condition of the same type on $f$.

We integrate over $X_{\infty}$ the function $F=f\circ F_Y$.
 Computing the integral 
from the definition, one can see that 
$\int_{X_{\infty}}e^{-F}=S_1-S_2$, where
$$S_1=\sum_{s\geq 0}E(Y_{s-1};u,v)(uv)^{-sn-f(s)},$$
$$S_2=\sum_{s\geq 0}E(Y_s;u,v)(uv)^{-(s+1)n	-f(s)}.$$

It is clear that every monomial which appears in the
$s$th term of $S_1$ has degree bounded above by 
$2P_1(s)$ and below by $2P_2(s)$, where
$P_1(s)=\dim\,Y_{s-1}-sn-f(s)$, and
$P_2(s) =-sn-f(s)$, for every $s\geq 0$ (recall our convention
that $Y_{-1}=X$). We have precisely one monomial of degree
$2P_1(s)$, namely $(uv)^{P_1(s)}$ whose coefficient is
$l(Y_{s-1})$, the number of irreducible components of
$Y_{s-1}$ of maximal dimension.

Similarly, every monomial which appears in the
$s$th term of $S_2$ has degree bounded above by $2Q_1(s)$ and
below by $2Q_2(s)$, where 
$$Q_1(s)=\dim\,Y_s-(s+1)n-f(s),$$
$$Q_2(s)=-(s+1)n-f(s).$$
There is exactly one monomial of degree $2Q_1(s)$, namely
$(uv)^{Q_1(s)}$ with coefficient $l(Y_s)$.

From the above evaluation of the terms in $S_1$
and $S_2$ and Lemma~\ref{measurable} we see that
$F$ is integrable. Moreover, it follows from condition
$(\star)$ that $P_1(s+1)<\min\{P_2(s), Q_2(s)\}$
for every $s$. We have also $Q_1(s)\leq P_1(s)$
for every $s$, with equality if and only if $s\geq 1$
and $\dim\,Y_s=\dim\,Y_{s-1}+n$. Lemma~\ref{measurable}
implies therefore that we have strict inequality
for infinitely many $s$. This shows that in $\int_{X_{\infty}}e^{-F}$
we have monomials of degree $2P_1(s)$ for infinitely many values of $s$.

We apply now the change of variable formula
in Proposition~\ref{change_of_variable} for the morphism $\pi$
to get
\begin{equation}
\int_{X_{\infty}}e^{-F}=\int_{X'_{\infty}}e^{-(F\circ\pi_{\infty}+F_W)}.
\end{equation}

Using the fact that $F\circ\pi_{\infty}=f\circ F_{\pi^{-1}(Y)}$ and
that $W\cup\pi^{-1}(Y)$ has simple normal crossings, one can compute
 explicitely
the integral. For a subset $J\subseteq\{1,\ldots,r\}$ let
$E_J^{\circ}=\cap_{i\in J}\setminus \cup_{i\not\in J}E_i$.
With this notation we have
$$\int_{X_{\infty}}e^{-F}=\sum_{J\subset\{1,\ldots,r\}}S_J,\,{\rm where}$$

$$ S_J=\sum_{\alpha_i\geq 1,i\in J}E(E_J^{\circ};u,v) (uv-1)^{|J|}
\cdot (uv)^{-n-\sum_{i\in J}\alpha_i(b_i+1)-f(\sum_{i\in J}\alpha_ia_i)}.$$

Every monomial in the term of $S_J$ corresponding to $(\alpha_i)_{i\in J}$
has degree bounded above by $2R_1(\alpha_i;i\in J)$ and below by
$2R_2(\alpha_i;i\in J)$, where
$$R_1(\alpha_i;i\in J)=-\sum_{i\in J}\alpha_i(b_i+1)-f(\sum_{i\in J}
\alpha_ia_i)$$
and $R_2(\alpha_i;i\in J)=R_1(\alpha_i;i\in J)-n$. Note that
 $R_1(\emptyset)=
-f(0)$.

Let us introduce the notation
$\tau(s)=\dim\,Y_s-(s+1)(n-q)$. We see that if $J\neq\emptyset$, then
$$R_1(\alpha_i; i\in J)= P_1(\sum_{i\in J}\alpha_ia_i)-
\tau(\sum_{i\in J}\alpha_i a_i-1) +\sum_{i\in J}\alpha_i(qa_i-b_i-1).$$
Moreover, property $(\star)$ implies that if $J\neq\emptyset$, then
\begin{equation}\label{eq1}
P_1(\sum_{i\in J}\alpha_ia_i+1)<R_2(\alpha_i;i\in J)
< R_1(\alpha_i;i\in J)<
\end{equation}
$$\min\{P_2(\sum_{i\in J}\alpha_ia_i-1), Q_2(\sum_{i\in J}\alpha_ia_i-1)\},$$
and $P_1(1)<R_2(\emptyset)$. This shows that the only monomial of the
form $(uv)^{P_1(s)}$ which could appear in the part of $S_J$ corresponding
to $(\alpha_j)_j$ is for $s=\sum_{j\in J}\alpha_ja_j$.

We prove now $(1)$.
 Suppose first that $(X, q\cdot Y)$ is log canonical,
 so that $b_i+1\geq qa_i$, for all $i$. Assume that for some $m\geq 0$,
we have $\tau(m)>0$. It follows from the above discussion that 
$(uv)^{P_1(m+1)}$ does not appear in $\int_{X_{\infty}}e^{-F}$.
As we have seen, this implies $\dim\,Y_{m+1}=\dim\,Y_m+n$. In particular, 
we have $\tau(m+1)>0$. Continuing in this way, we deduce
$\dim\,Y_{s+1}=\dim\,Y_s+n$ for all $s\geq m$, which is impossible.

Conversely, suppose that for some fixed $m$
such that $a_i\mid (m+1)$ for all $i$ 
we have $\tau(m)\leq 0$, and that 
$b_j+1<ca_j$, for some $j$.

We choose the function $f$ such that in addition to $(\star)$ it satisfies
the following condition. For every $p$, let ${\mathcal J}_p$ be the finite
set consisting of all the pairs $(J, (\alpha_i)_{i\in J})$ such that
$\sum_{i\in J}a_i\alpha_i=p$. The extra condition we require for $f$
is that for every $(J, (\alpha_i)_{i\in J})\in {\mathcal J}_{m+1}$
and every $(J', (\alpha'_i)_{i\in J'})\in {\mathcal J}_p$, with $p\leq m$,
we have
\begin{equation}\label{requirement}
f(m+1)> f(p)-\sum_{i\in J}\alpha_i(b_i+1)+\sum_{i\in J'}\alpha'_i(b_i+1)+n.
\end{equation}
This means that $R_1(\alpha_i;i\in J)<R_2(\alpha'_i;i\in J')$.

Note that from $(\ref{eq1})$ we can deduce that if 
$(J',(\alpha'_i)_{i\in J'})\in {\mathcal J}_p$, with $p\geq m+2$, then
\begin{equation}\label{eq2}
P_1(m+1)> R_1(\alpha'_i;i\in J').
\end{equation}

On the other hand, the top degree monomials
which appear in different terms of the sums $S_J$ (for possible
different $J$) don't cancel each other, as they have positive coefficients.
Let $d$ be the highest degree of a monomial which appears in a term
corresponding to some $(J, (\alpha_i)_{i\in J})\in {\mathcal J}_{m+1}$.
The previous remark shows that the corresponding monomial does not
cancel with the other monomials which appear in terms corresponding to 
 $(J', (\alpha'_i)_{i\in J'})\in {\mathcal J}_{m+1}$. 

Obviously we have $(\{j\}, (m+1)/a_j)\in {\mathcal J}_{m+1}$ and
therefore our assumption that $\tau(m)\leq 0$ and $b_j+1<qa_j$
gives $d>2P_1(m+1)$. We also deduce from $(\ref{eq1})$ that 
$d<\min\{2P_2(m), 2Q_2(m)\}$. 

We see from $(\ref{requirement})$ and $(\ref{eq2})$ that this monomial of
 degree $d$ does not cancel with other monomials in terms corresponding to
$(J', (\alpha'_i)_{i\in J'})\in {\mathcal J}_p$, if $p\neq m+1$.
Therefore in the integral
 $\int_{X_{\infty}}e^{-F}$ we have a monomial of degree $d$, with
$2P_1(m+1)<d<\min\{2P_2(m), 2Q_2(m)\}$, a contradiction.

The proof of $(2)$ is entirely similar.
\end{proof}

In order to state the formula for the log canonical threshold which
follows from Theorem~\ref{main}, we make the following definition.
Consider a pair $(X,Y)$ as before and $Z\subseteq X$
 a nonempty closed subset.
We allow also the case $Y=X$.

\begin{definition}
With the above notation,
if $m\in\NN$, we define
$\dim_ZY_m$ to be the dimension of $Y_m$ along $Y_m\cap\rho_m^{-1}(Z)$ i.e.
the maximum dimension of an irreducible component $T$ of $Y_m$
such that $\rho_m(T)\cap Z\neq\emptyset$.
\end{definition}

\begin{remark}\label{existence_U}
It follows from our discussion of jet schemes in the previous section that
 for every irreducible component
$T$ of $Y_m$, the image $\rho_m(T)$ is closed in $Y$ (and therefore in $X$).
Indeed, since $T$ must be invariant 
with respect to the ``action'' of ${\bf A}^1$
on $Y_m$, it follows that $\rho_m(T)=(\sigma_m^Y)^{-1}(T)$.

 Therefore there is an
open neighbourhood $U$ of $Z$ such that for
every irreducible component $T$ of $Y_m$ with $\rho_m(T)\cap Z=\emptyset$,
$U$ does not intersect $\rho_m(T)$.
For every such $U$, we have $\dim_ZY_m=\dim\,(Y\cap U)_m$.
\end{remark}

\smallskip

In the corollary below, we keep the previous notation for
a log resolution of the pair $(X,Y)$.

\begin{corollary}\label{compute_lc}
For every pair $(X,Y)$ and $Z\subseteq X$ as above,
$$c_Z(X,Y)=n-\sup_{m\geq 0}{{\dim_ZY_m}\over {m+1}}.$$ Moreover, if
$m\geq 0$ is such that $a_i\mid(m+1)$, for all $i$
such that $\pi(E_i)\cap Z\neq\emptyset$, then
$$c_Z(X,Y)=n-{{\dim_ZY_m}\over {m+1}}.$$
\end{corollary}

\begin{proof}
Note first that the formula is trivial when $Y=X$, since in this case
for every $m\geq 1$, we know that $Y_m=X_m$
is a smooth variety  of dimension $(m+1)\,\dim\,X$.
Suppose from now on that $Y\neq X$.

For a fixed $m$, we can find an open neighbourhood $U$ of $Z$
such that $\dim_ZY_m=\dim\,(Y\cap U)_m$ and $c_Z(X,Y)=c(U,Y\cap U)$.
Therefore, in order to prove the second assertion of the corollary we
may assume $Z=X$, in which case the formula follows from Theorem~\ref{main}.

In order to prove the first assertion, it remains to be proved
that $c_Z(X,Y)\geq n-(\dim_ZY_m)/(m+1)$, for every $m$. We reduce again
to the case when $Z=X$, when we conclude by applying Theorem~\ref{main}.
\end{proof}

\begin{corollary}\label{independance}
If $(X,Y)$ is a pair and $Z$ a nonempty closed subset of $X$, then $c_Z(X,Y)$
depends only on $Y$, $\dim\,X$ and the closed subset $Y\cap Z\subseteq Y$,
but not on the particular ambient variety $X$ or the embedding 
$Y\hookrightarrow X$.
\end{corollary}

\begin{proof}
The statement is a consequence of the formula in Corollary~\ref{compute_lc}.
\end{proof}

\smallskip

We give now a version of the formula in Corollary~\ref{compute_lc} involving
only the fibers of $\rho_m^Y$
 over $Y\cap Z$. The drawback is that in this case,
the supremum which appears in the formula can not be obtained, in general,
for a specific integer.

\begin{corollary}\label{compute_lc2}
If $(X,Y)$ is a pair and $Z\subset X$ a closed subset as before, then
$$c_Z(X,Y)=\dim\,X-\sup_{m\geq 0}
{{\dim\,(\rho^Y_m)^{-1}(Y\cap Z)}\over{m+1}}.$$
\end{corollary}

\begin{proof}
Since the cases $Y=X$ and $Y\cap Z=\emptyset$ are trivial, we may assume 
that we are in neither of them.
We obviously have $\dim\,(\rho^Y_m)^{-1}(Z\cap Y)\leq\dim_ZY_m$,
for every $m$,  so that
the formula in Corollary~\ref{compute_lc} gives
$$c_Z(X,Y)\leq\dim\,X-\sup_{m\geq 0}
{{\dim\,(\rho^Y_m)^{-1}(Z\cap Y)}\over{m+1}}.$$

On the other hand, the formula giving the lower bound for the dimension
of the fibers of a
 morphism implies that 
$$\dim_ZY_m\leq\dim\,Y+\dim\,(\rho^Y_m)^{-1}(Y\cap Z).$$

Note that by Corollary~\ref{compute_lc}, there is
a positive integer $N$ such that $c_Z(X,Y)=\dim\,X-{{\dim_ZY_m}\over
{m+1}}$ if $N\,|(m+1)$. Therefore we deduce
$$c_Z(X,Y)\geq\dim\,X-{{\dim\,Y}\over{m+1}}-
{{\dim\,(\rho^Y_m)^{-1}(Y\cap Z)}\over{m+1}},$$
 for every $m$, such that $N\,|(m+1)$.
This implies that for every $p\geq 1$,
$$c_Z(X,Y)\geq \dim\,X- {{\dim\,Y}\over {pN}}-\sup_{m\geq 0}
{{\dim\,(\rho^Y_m)^{-1}(Y\cap Z)}\over{m+1}},$$ so that in fact
$$c_Z(X,Y)\geq \dim\,X-\sup_{m\geq 0}
{{\dim\,(\rho_m^Y)^{-1}(Y\cap Z)}\over{m+1}},$$
which completes the proof of the corollary.
\end{proof}

\smallskip

We end this section with the following
\begin{example}
Consider the case of a cusp:
$$Y=Z(u^2-v^3)\hookrightarrow X=\AA^2.$$
$Y$ is a curve and $Y_{\rm sing}=\{0\}$, so that $\dim\,Y_m=\dim(\rho_m^Y)
^{-1}(0)$.

 In order to describe $(\rho_m^Y)^{-1}(0)$, note that it
consists of ring homomorphisms
$$\phi\,:\,\CC[u,v]/(u^2-v^3)\longrightarrow \CC[t]/(t^{m+1})$$
such that ${\rm ord}(\phi(u))$, ${\rm ord}(\phi(v))\geq 1$. If $m>5$,
this implies that $\phi(u)=t^3f$ and $\phi(v)=t^2g$, such that the classes
$\overline{f}$, $\overline{g}$ of $f$ and $g$ in 
$\CC[t]/(t^{m-5})$ satisfy $\overline{f}^2-\overline{g}^3=0$.
In this way we get an isomorphism $(\rho_m^Y)^{-1}(0)\simeq
Y_{m-6}\times\AA^7$, so that $\dim\,Y_m=\dim\,Y_{m-6}+7$.
This easily gives ${\rm sup}_m\{(\dim\,Y_m)/(m+1)\}=7/6$, so that
$c(\AA^2,Y)=5/6$.
\end{example}

\bigskip
\bigskip

\section{The log canonical threshold via jets}

\bigskip

In this section we apply Corollary~\ref{compute_lc} to deduce properties
of the log canonical threshold which appear in the analytic context in 
\cite{demailly}. In the codimension one case, some of these properties
are proved also in \cite{kollar} using log resolutions. We believe that our
treatment, using properties of the jet schemes, simplifies many of the proofs.

In this section, when we consider pairs $(X,Y)$ we allow also the case $Y=X$.
Recall that we follow the convention that $\dim(\emptyset)=-\infty$.

\begin{proposition}(\,\cite{demailly} 1.4)\label{inclusion}
If $X$ is a smooth variety, $Y'$, $Y''$ are two closed subschemes
such that $Y'\hookrightarrow Y''$ 
and $Z\subseteq X$ is a nonempty closed subset, then $c_Z(X,Y')\geq
c_Z(X, Y'')$.
\end{proposition}

\begin{proof}
Since $Y'$ is a subscheme of $Y''$, it follows that $Y'_m$ is a subscheme
of $Y''_m$, so that $\dim_ZY'_m\leq\dim_ZY''_m$ for every $m$.
 The assertion now follows
from Corollary~\ref{compute_lc}.
\end{proof}

\begin{proposition}(\,\cite{demailly} 1.4)\label{bound1}
For every pair $(X,Y)$
 and every nonempty closed subset
$Z\subseteq X$, we have $c_Z(X,Y)\leq {\rm codim}_Z(Y,X)$,
where ${\rm codim}_Z(X,Y)$ denotes the smallest codimension of an irreducible
component of $Y$ meeting $Z$.
\end{proposition}

\begin{proof}

Using Corollary~\ref{compute_lc}, we have
$$c_Z(X,Y)=\dim\,X-\sup_{m\geq 0}{{\dim_ZY_m}\over{m+1}}$$
$$\leq\dim\,X-\dim_ZY_0={\rm codim}_Z(Y,X).$$
\end{proof}

For a pair $(X,Y)$ and a point $x\in X$, we define
${\rm mult}_xY$ as follows. If $\cO_{X,x}$ is the local ring
of $X$ at $x$, with maximal ideal $\underline{m}_{X,x}$,
and if ${\mathcal I}_{Y,x}\subseteq\cO_{X,x}$ is the ideal of $Y$ at $x$,
then ${\rm mult}_xY$ is the largest $q\in\NN$ such that
$I_{Y,x}\subseteq\underline{m}_{X,x}^q$. Note that if
$Y$ is a divisor, then this is the same with the Hilbert-Samuel
multiplicity, but in general the two notions are different.

\begin{proposition}(\,\cite{demailly} 1.4, \cite{kollar} 8.10)\label{bound2}
Let $(X,Y)$ be a pair and $Z\subseteq X$ a nonempty closed subset.
If $q$ is the smallest $p\in\NN$ such that 
${\rm mult}_xY\leq p$ for every $x\in Z$, then
\begin{equation}
1/q\leq c_Z(X,Y)\leq (\dim\,X)/q.
\end{equation}
\end{proposition}

\begin{proof}
The cases $Y=X$ (when $q=\infty$) and $Y\cap Z=\emptyset$
(when $q=0$) are trivial, so that we may suppose that
we are in neither of them.

There is $y\in Z$ such that ${\rm mult}_{y}Y=q$. If
$\rho_{q-1}\,:\,X_{q-1}\longrightarrow X$ is the canonical projection,
then $\rho_{q-1}^{-1}(y)\subseteq Y_{q-1}$. Indeed, every
local morphism
$\cO_{X,y}\longrightarrow\CC[t]/(t^q)$ factors through $I_{Y,y}$,
as $I_{Y,y}\subseteq\underline{m}_{X,y}^q$.

This gives $\dim_ZY_{q-1}\geq\dim\,\rho_{q-1}^{-1}(y)=(q-1)\dim\,X$.
Corollary~\ref{compute_lc} implies that
$c_Z(X,Y)\leq\dim\,X-(1/q)(q-1)\dim\,X=(\dim\,X)/q$.

For the other inequality we use Lemma~\ref{measurable}. 
This shows that for every $x\in X$ and every $m\in\NN$, we have
$\dim(\rho_m^{-1}(x)\cap Y_m)\leq m\dim\,X-[m/{\rm mult}_xY]$. Since
${\rm mult}_xY\leq q$ for every $x\in Z$, the same is true on an
open neighbourhood $U$ of $Z$. We deduce
$$\dim\,(U\cap Y)_m\leq \dim\,Y+m\cdot\dim\,X-[m/q], $$
for every $m$, 
and an easy computation gives $\dim\,(Y\cap U)_m/(m+1)\leq
\dim\,X-1/q$. Corollary 
~\ref{compute_lc} implies $c_Z(X,Y)\geq c(U, U\cap Y)
\geq 1/q$.
\end{proof}

\begin{proposition}(\,\cite{demailly} 2.7)\label{product}
If $(X',Y')$ and $(X'',Y'')$ are two pairs and $Z'\subseteq X'$,
$Z''\subseteq X''$ are two nonempty closed subsets, then
\begin{equation}
c_{Z'\times Z''}(X'\times X'', Y'\times Y'')=
c_{Z'}(X',Y')+c_{Z''}(X'',Y'').
\end{equation}
\end{proposition}

\begin{proof}
Proposition~\ref{product_of_jets} gives a canonical isomorphism
$(Y'\times Y'')_m\simeq Y'_m\times Y''_m$, for every $m$. 
Therefore we have $\dim_{Z'\times Z''}(Y'\times Y'')_m=
\dim_{Z'}Y'_m+\dim_{Z''}Y''_m$. We pick $m$ such that
$a'_i$, $a''_j$, $a_k\mid (m+1)$
for all $i$, $j$ and $k$, where $a'_i$, $a''_j$
and $a_k$
are the coefficients appearing in log resolutions of $(X',Y')$,
$(X'',Y'')$ and $(X'\times X'', Y'\times Y'')$.
 We conclude by applying the second assertion
in Corollary~\ref{compute_lc}. 
\end{proof}

\begin{proposition}(\,\cite{demailly} 2.2)\label{inverse_adjunction}
Let $(X,Y)$ be a pair. 
If $H\subset X$ is a smooth irreducible divisor 
 and $Z\subseteq H$ is a nonempty closed subset, then
\begin{equation}
c_Z(X,Y)\geq c_Z(H, Y\cap H).
\end{equation}
\end{proposition}

\begin{proof}
Let $U$ be an open neighbourhood of $Z$ in $X$ such that
$c_Z(H, H\cap Y)=c\,(U\cap H, U\cap H\cap Y)$.
Since $c_Z(X,Y)\geq c_{U\cap H}
(U, U\cap Y)$, by restricting everything to
$U$, we may assume that $Z=H$.

From the local description of the jet schemes, it follows that since
$H\cap Y$ is defined locally in $Y$ be one equation, we have
$(H\cap Y)_m$ defined locally in $Y_m$ by $m+1$ equations, 
for every $m\in\NN$. Moreover, if $\rho_m\,:\,X_m\longrightarrow X$
is the canonical projection, then for every
irreducible component $T$ of $Y_m$ such
that $\rho_m(T)\cap H\neq\emptyset$, we have $T\cap (Y\cap H)_m\neq
\emptyset$.

Indeed, it follows from our discussion in Section 2 that if
$\sigma_m^Y\,:\,Y\longrightarrow Y_m$ is the ``zero section'',
then
$\sigma^Y_m(\rho_m(T))\subseteq T$.
Moreover, since  the ``zero section'' is functorial, we have 
$\sigma^Y_m(Y\cap H)\subset (Y\cap H)_m$, and we deduce
that $(Y\cap H)_m\cap T$ contains $\sigma^Y_m(\rho_m(T)\cap H)$.

It follows that $\dim_HY_m\leq\dim\,(Y\cap H)+m+1$, so 
that
 $$\sup_{m\geq 0}{{\dim_HY_m}\over {m+1}}
\leq\sup_{m\geq 0}{{\dim\,(Y\cap H)_m}\over {m+1}}+1.$$
Applying Corollary~\ref{compute_lc}, we deduce the inequality in the
proposition.
\end{proof}

\begin{corollary}\label{gen_codimension}
The statement of Proposition~\ref{inverse_adjunction}
remains true if we replace $H$ with a smooth variety of arbitrary 
codimension.
\end{corollary}

\begin{proof}
If we cover $X$ with open subsets $U_i$, then we obviously have
$c_Z(X,Y)=\inf_i\{c_{Z\cap U_i}(U_i, Y\cap U_i)\}$ and
a similar relation for the restrictions to $H$. Since $H$ is smooth,
it is locally a complete intersection, so that we can find an open cover of
$X$ such that on each subset $U_i$, $H\cap U_i$ is an intersection of
smooth divisors. It is therefore enough to apply Proposition~\ref
{inverse_adjunction} inductively on each of these open subsets.
\end{proof}

This corollary can be used to deduce a formula for the log canonical
threshold of an intersection $Y'\cap Y''$ in terms of the dimensions
of the jet schemes of $Y'$ and $Y''$. We give below the case when
the closed subset $Z\subseteq X$ is a point, when the formula
has a particularly nice form.

\begin{proposition}(\,\cite{demailly} 2.9)\label{intersection}
If $X$ is a smooth variety, $Y'$, $Y''$
two proper closed subschemes and $z\in Y'\cap Y''$ a point, then:
\begin{equation}
c_z(X,Y')+c_z(X,Y'')\geq c_z(X,Y'\cap Y'').
\end{equation} 	
\end{proposition}

\begin{proof}

Consider the following cartesian diagram:
\[
\begin{CD}
Y'\cap Y'' @>>> X \\
@VVV @VV{\Delta}V \\
Y'\times Y'' @>>>X\times X
\end{CD}
\]
where the horizontal maps
 are the natural embeddings and $\Delta$ is the
diagonal embedding. Applying Corollary~\ref{gen_codimension} to the
embedding $\Delta$, the pair $(X\times X, Y'\times Y'')$,
and the closed subset $\{(z,z)\}\subset \Delta(X)$, we deduce
$$c_{(z,z)}(X\times X, Y'\times Y'')\geq c_z(X, Y'\cap Y'').$$

On the other hand, Proposition~\ref{product} implies
$$c_{(z,z)}(X\times X,Y'\times Y'')=c_z(X,Y')+c_z(X,Y''),$$
and these two relations prove the proposition.
\end{proof}

We use now Corollary~\ref{compute_lc} and Proposition~\ref{semicontinuity1}
to prove a result of Demailly and Koll\'{a}r (\,\cite{demailly}) about the
semicontinuity of the log canonical threshold. We first need 
a preliminary result. For the notation concerning the coefficients
of the divisors in a log resolution, see the beginning of the previous section.

\begin{lemma}\label{finite}
Let $\pi\,:\,{\mathcal X}\longrightarrow S$ be a smooth morphism and
${\mathcal Y}\hookrightarrow {\mathcal X}$ a closed subscheme. It is
possible to find log resolutions for each fiber $({\mathcal X}_s,
{\mathcal Y}_s)$ such that the coeficients $a_i$, $b_i$ which appear in all
these resolutions form a finite set.
\end{lemma}

\begin{proof}
We make induction on $\dim\,S$. It is enough to find a locally
closed cover of $S$ such that the restriction of $\pi$ over each member of the
cover has the desired property. Therefore we may assume that $S$ is
a smooth variety and it is enough to find a nonempty open subset $U$ of $S$
over which the restriction of $\pi$ has this property.
 By a theorem of Nagata, we can embed ${\mathcal X}$
as an open subscheme in a scheme ${\mathcal X}'$ proper over $S$.
By replacing ${\mathcal X}$ with ${\mathcal X}'$ and ${\mathcal Y}$
by its closure in ${mathcal X}'$, we may assume that $\pi$ is proper.
After restricting to an open subset of $S$, we may assume that $\pi$
is also smooth.

 Let $\pi'\,:\,
\widetilde{\mathcal X}\longrightarrow {\mathcal X}$ be a log resolution for 
the pair $({\mathcal X}, {\mathcal Y})$. It is easy to see that after 
further restricting over an open subset of $S$ we can assume that for
every $s\in S$, the restriction to the fiber $\pi_s\,:\,
\widetilde{{\mathcal X}}
_s\longrightarrow {\mathcal X}_s$ is a log resolution of $({\mathcal X}_s,
{\mathcal Y}_s)$, in which case the assertion is obvious.
\end{proof}

We can give now the proof of the semicontinuity result.

\begin{theorem}(\,\cite{demailly} 3.1)\label{semicontinuity2}
Let $\pi\,:\,{\mathcal X}\longrightarrow S$ be a smooth morphism,
${\mathcal Y}\hookrightarrow {\mathcal X}$ a closed subscheme, and
$\tau\,:\,S\longrightarrow {\mathcal Y}$ a section of $\pi|_{\mathcal Y}$.
The function defined by
$$f(s)=c_{\tau(s)}({\mathcal X}_s,{\mathcal Y}_s),$$
for every closed point $s\in S$, is lower semicontinuous.
\end{theorem}

\begin{proof}
Lemma~\ref{finite} and Proposition~\ref{lc_formula} show
that the set 
 $\{c_{\tau(s)}({\mathcal X}, {\mathcal Y})\,|\,s\in S\}$
is finite. Moreover, it follows from Lemma~\ref{finite} and
Corollary~\ref{compute_lc} that we can find $N\geq 1$ such that
for every $s\in S$ and every $m$ with $N| (m+1)$, we have
$$c_{\tau(s)}({\mathcal X}_s, {\mathcal Y}_s)=
\dim\,{\mathcal X}_s-{{\dim_{\tau(s)}({\mathcal Y}_s)_m}\over{m+1}}.$$

Fix $s_0\in S$. Since there are only finitely many log canonical
thresholds to consider,
 it is enough to show that for every $\epsilon>0$, there is
an open neighbourhood $U$ of $s_0$ such that $c_{\tau(s)}({\mathcal X}_s,
{\mathcal Y}_s)\geq c_{\tau(s_0)}({\mathcal X}_{s_0},{\mathcal Y}_{s_0})
-\epsilon$,
for every $s\in U$. Since $\pi$ is flat, by restricting to an
open neighbourhood of $s_0$ we may assume that $\dim\,{\mathcal X}_s$
is constant on $S$. 

Therefore it is enough to find $m$ with $N|\,(m+1)$ and $U$ such that
$$\dim_{\tau(s)}({\mathcal Y}_s)_m/(m+1)\leq \dim
_{\tau(s_0)}({\mathcal Y}_{s_0})_m/(m+1)+\epsilon,$$
for all $s\in U$.
We fix $m$ such that $N|\,(m+1)$ and $\dim\,{\mathcal X}_s/(m+1)
\leq\epsilon$ for all $s$. Using Proposition~\ref{semicontinuity1},
we choose an open neighbourhood $U$ of $s_0$ such that
$$\dim\,(\rho_m^{{\mathcal Y}_s})^{-1}(\tau(s))\leq \dim\,(\rho
^{{\mathcal Y}_{s_0}}_m)^{-1}(\tau(s_0)),$$
for all $s\in U$. The following inequalities show that $U$ satisfies
our requirement:
$$\dim_{\tau(s)}({\mathcal Y}_s)_m/(m+1)\leq
\dim\,(\rho_m^{{\mathcal Y}_s})^{-1}(\tau(s))/(m+1)+\dim\,{\mathcal X}_s/
(m+1)\leq$$
$$\dim\,(\rho_m^{{\mathcal Y}_{s_0}})^{-1}(\tau(s_0))/(m+1)+\epsilon
\leq \dim_{\tau(s_0)}({\mathcal Y}_{s_0})_m/(m+1)+\epsilon.$$
\end{proof}

\smallskip

We conclude by showing how the formula in Corollary~\ref{compute_lc}
can be used to explicitely compute the log canonical threshold. We consider
the case of monomial ideals and derive a formula from \cite{agv}
(see also \cite{howald}
for a formula in this context for all the multiplier ideals).

Let us fix the notation. $X={\bf A}^n$ is the affine space and $Y=V(I)$,
where $I\subset R=k[X_1,\ldots, X_n]$ is a monomial ideal such that
$(0)\neq I\neq R$. For ${\bf a}\in\ZZ^n$, we use the notation
$X^{\bf a}=\prod_iX_i^{a_i}$. The vector $(1,\ldots,1)\in\ZZ^n$
is denoted by ${\bf e}$.
 The Newton polytope of $I$, denoted by $P_I$, is the
convex hull (in $\RR^n$) of the set $\{{\bf a}\in\ZZ^n\,|\,
X^{\bf a}\in I\}$. 

\begin{proposition}(\cite{agv}, \cite{howald})\label{monomial}
With the above notation, the log canonical threshold of the pair $(X,Y)$
is given by
$$c(X,Y)=\sup\{r>0\,|\,{\bf e}\in rP_I\}.$$
\end{proposition}

\begin{proof}
We will use the formula in Corollary~\ref{compute_lc}, so that we estimate
first $\dim\,Y_m$, for every $m\geq 1$. Note that $Y_m$ can be covered
by the locally closed subsets $Z_{a_1,\ldots,a_n}$, with
$0\leq a_i\leq m+1$ for all $i$, where $Z_{a_1,\ldots,a_n}$
is the set of ring homomorphisms
$\phi\,:\,R/I\longrightarrow k[t]/(t^{m+1})$, with ${\rm ord}\,(\phi(X_i))
=a_i$, for every $i$. We put ${\rm ord}(\phi(X_i))=m+1$
if $\phi(X_i)=0$.

It is clear that if $Z_{a_1,\ldots,a_n}\neq\emptyset$, then
$\dim\,Z_{a_1,\ldots,a_n}=(m+1)n-\sum_i a_i$. On the other hand,
$Z_{a_1,\ldots,a_n}\neq\emptyset$ if and only if for every
${\bf b}\in\NN^n$ such that $X^{\bf b}\in I$, we have $\sum_ia_ib_i
\geq (m+1)$.

Let $P_I^{\circ}$ be the polar polyhedron of $P_I$, defined by
$$P_I^{\circ}=\{{\bf u}\in\RR^n\,|\,\sum_iu_iv_i\geq 1, {\rm for}\,{\rm all}\,
{\bf v}\in P_I\}.$$ 
We see that $Z_{a_1,\ldots,a_n}\neq\emptyset$ if and only if
$(a_1/(m+1),\ldots, a_n/(m+1))\in P_I^{\circ}$.

From the above discussion, we deduce that
$${{\dim\,X_m}\over{m+1}}=n-{1\over{(m+1)}}\inf_{\bf a}\sum_ia_i,$$
where the infimum is taken over all ${\bf a}\in\NN^n\cap (m+1)P_I
^{\circ}$ such that $a_i\leq m+1$ for all $i$.
 The formula in Corollary~\ref{compute_lc} gives
$$c\,(X,Y)=\inf_{\bf a}\sum_ia_i,$$
the infimum being taken over all ${\bf a}\in P_I^{\circ}\cap\QQ^n\cap 
[0,1]^n$.

We clearly have $P_I^{\circ}\subset\RR_+^n$.
Note that if ${\bf a}\in P_I^{\circ}$ and ${\bf a}'$ is defined
by $a'_i={\rm min}\{a_i,1\}$ for all $i$, then ${\bf a}'\in
P_I^{\circ}$. We deduce that
 $c\,(X,Y)=\inf_{{\bf a}\in P_I^{\circ}}\sum_ia_i$. In order to
complete the proof it is enough to note that ${\bf e}\in rP_I$
if and only if for every ${\bf a}\in P_I^{\circ}$ we have $\sum_i
a_i\geq r$ (this comes from the fact that $P_I$ is the polar
polyhedron of $P_I^{\circ}$.
\end{proof}

\bigskip


\providecommand{\bysame}{\leavevmode\hbox to3em{\hrulefill}\thinspace}

\end{document}